\documentclass[12pt]{article}
\usepackage{amsmath,amscd,amsthm,amssymb}

 \newtheorem{theorem}{Theorem}[section]
 
 \newtheorem{lemma}{Lemma}[section]
 
 \theoremstyle{definition}
 
 \theoremstyle{remark}

 \numberwithin{equation}{section}
\begin {document}

\begin{center}
\bf{ON SELFADJOINT EXTENSIONS OF SYMMETRIC OPERATOR WITH EXIT FROM SPECE}
\end{center}
\begin{center}
Khalig M.Aslanov$^{a}$, Nigar M. Aslanova$^{b,c}$
			\\
$^{a}$\it Azerbaujan State Economic University (UNEC),\\
6, Istiglaliyyat str., AZ 1001, Baku, Azerbaijan.\\

$^{b,c}$\it Institute of Mathematics and Mechanics,\\ National Academy of Sciences of Azerbaijan Az-1141, Baku, Azerbaijan, \\
Azerbaijan University of Architecture and Construction\\

{\it E-mail addresses}: xaliqaslanov@yandex.ru, nigar.aslanova@yahoo.com

\end{center}

\begin{abstract}
We investigate minimal operator corresponding to operator differential expression with exit from space, study its selfajoint extensions, also for one particular selfadjoint extension, corresponding to boundary value problem with some rational function of eigenparameter in boundary condition establish asymptotics of spectrum and derive trace formula.
\end{abstract}

\textbf{2010 Mathematics Subject Classifications:} 34B05, 34G20, 34L20, 34L05, 47A05, 47A10.

\textbf{Keywords: }Hilbert space, differential operator equation, selfadjoint extensions with exit from space, spectrum, eigenvalues, trace class operators, regularized trace.

\section{Introduction}

Let $H$ be separable Hilbert space and $L_{2} \left(H,\left(0,1\right)\right)$ set of vector functions $y\left(t\right)\ $with values from $H$for each $t\in \left[0,1\right]$ satisfying $\int _{0}^{1}\left\| y\left(t\right)\right\| ^{2} dt<\infty  $($\left\| \cdot \right\| $-z is norm in $H$). $L_{2} \left(H,\left(0,1\right)\right)$ is Hilbert space with scalar product defined as
\[\left(y\left(z\right),z\left(t\right)\right)_{L^{2} \left(H,\left(0,1\right)\right)} =\int _{0}^{1}\left(y\left(t\right),z\left(t\right)\right) dt,\]
$\left(\cdot ,\cdot \right)$ is a scalar product in $H$). We at  first will consider  operator corresponding to differential expression $ly\equiv -y''+Ay+q\left(t\right)y$ in the ``larger space'', namely, direct sum Hilbert space   $L_{2} =L_{2} \left(H,\left(0,1\right)\right)\oplus H^{2} $ elements of which are triples $Y=\left\{y\left(t\right),y_{1} ,y_{2} \right\},$ $y\left(t\right)\in L_{2} \left(H,\left(0,1\right)\right),$ $y_{1} ,y_{2} \in H$ . About coefficients of  $l$ we assume that $A=A^{*} >I$, $A^{-1} \in \sigma _{\infty } $($\sigma _{\infty } $ class of compact operators), $I$- is identity operator in $H$. Operator valued function $q(t)$ is selfadjoint, bounded operator in $H$ and $\left\| q(t)\right\| \le const$. From that it follows $q(t)$ is bounded also in $L_{2} \left(H,(0,1)\right)$.

Let M=$D\left(L_{0}^{\prime}
\right)=\left\{Y:Y=\left\{y\left(t\right),y_{1},y\left(1
\right)\right\},y(t)\in C_{0}^{\infty } \left(H_{\infty },
\left[0,1\right]\right)\right\},$$ $$H_{\infty } =\bigcap _{j=1}^{\infty }D\left(A^{j} \right) $,  $D\left(A^{j}
\right)$ denotes domain of definition of $A^{j} $,
$C_{0}^{\infty } \left(H_{\infty},
\left[0,1\right]\right)$ is set of functions whose values
are from $H_{\infty } $, which  vanish in the vicinity of
zero and infinitely many time differentiable. It might
by be easily shown that M  is dense in $L_2.$ Define in
M=$D\left(L^{\prime}_{0} \right)$ operator $L^{\prime}_{0} $ by $L^{\prime}_{0}
Y=\left\{l\left(y\right),-y\left(1\right),y^{\prime}\left(1\right
)-y_{1} \right\}.$

Denote closure of $L'_{0} $ in $L_{2} $ by $L_{0} $ and call it minimal operator. Adjoint of $L_{0} $ call maximal operator.

We will show symmetricity of $L'_{0} $, study its adjoint operator which is called maximal operator, give description of its domain,  also  define selfadjoint extensions with discrete spectrum. For particular selfadjoint extension corresponding to boundary value problem  \eqref{eq-1_1}, \eqref{eq-1_2}

\begin{equation} \label{eq-1_1}
-y^{''} (x)+Ay(x)+q(x)y(x)=\lambda y(x),
\end{equation}

\begin{equation} \label{eq-1_2}
y(0)=0,{\kern 1pt} {\kern 1pt} {\kern 1pt} y^{'} (1)=(\lambda -\frac{1}{\lambda } )y(1).
\end{equation}
with rational function of spectral parameter  in boundary condition we obtain asymptotic formula for  distribution function of eigenvalues $N(\lambda )$ ($N(\lambda )=\sum _{\lambda _{n} <\lambda } 1$, here $\lambda \, _{n} $ are eigenvalues of operator associated with that problem), find also asymptotics of $\lambda _{n} $, and finally derive first order regularized trace formula. Some additional conditions on $q(x)$ needed for derivation of trace formula will be given later.

As it is seen from \eqref{eq-1_1},\eqref{eq-1_2} eigenparameter $\lambda $ appear in equation linearly, but boundary condition at 1 contains some Herglotz-Nevanlinna type rational function  of $\lambda $.

Note that in many works where considered spectral problems with eigen parameter in equation as well as in boundary condition, references to which we give below authors go in next  way. For concrete   spectral problem with eigenparameter included boundary condition they firstly linearize it by introducing larger space (direct sum space), so exit from original space, define operator in direct sum space corresponding to that problem , then show selfadjointness of that  operator and study other spectral questions of their interest. Our organization of paper depicted above is more closely related to \cite{Gor75} where abstract Sturm-Liouville problem is considered.

It is well known that second order ordinary differential equation involving eigenparameter first was studied by Liouville in 1837. At the beginning of last century G.D. Birkhoff \cite{Bir08} and J.D. Tamarkin \cite{Tam28} treated $n$-th order differential equation and system of first order differential equations including $n$-th powers of eigenparameter, where asymptotic behavior of solutions, expansion of some classes of functions in series of eigen and associated functions and other spectral questions were treated. After appeared lots of works in that direction. In some works from that lot different linearization techniques were used and extending space  by introducing new direct sum space    problem was reduced to spectral problem $Ay=\lambda y$ for some operator $A$ in that "larger spaces", where authors studied different spectral questions such as minimality, completeness or basicity of eigen or associate functions and so forth. Note here that linearization techniques were used for example, in works by Gohberg, Kaashoek, Lay \cite{Goh}, Shkalikov\cite{Shk}, Tretter \cite{Tre,Tre1} and many others. But we will distinguish through them those studies where after linearization problem was reduced to spectral problem for selfadjoint operator with exit from space, where equation is considered, because they are closely related to problem of our concern. In problems for Sturm-Liouville equations in scalar case with eigenvalue parameter involved in boundary condition linearly, polynomially or rationally authors indicate as original works the ones by Walter \cite{Wal} (in1973) and Fulton \cite{Ful}(in 1976). Thus, in that scalar cases problems are recast as $Ay=\lambda y$ in space
\begin{equation} \label{eq-1_3}
L_{2} =L_{2} (\Omega ){\mathop{\oplus }\limits_{}} C^{d}
\end{equation}
where $\Omega $ is interval where equation is satisfied, $d$ is number determined by degree of $\lambda $ in boundary condition when $\lambda $ polynomially involved in boundary condition, by number of $\lambda $ dependent boundary conditions and number of $\lambda $ dependent partial fractions and positive degrees of $\lambda $ if it involved rationally in boundary condition. $A$ becomes selfadjoint in Hilbert space $L_{2} $ or in Pontryagin space with respect to indefinite metrics $[\cdot ,\cdot ]=(\cdot ,\cdot )_{L_{2} (\Omega )} -(\cdot ,\cdot )_{C^{d} } $ introduced in $L_{2} $.

In \cite{Bin,Bin1} considered two point Strum-Liouville problem in scalar case with rational function of $\lambda $ in boundary condition at one end point of interval where equation is satisfied. Authors give nesessary and sufficinet condition for rational function with simple poles to be Herglotz-Nevanlinna function. In first article \cite{Bin} existence and asymptotics of eigenvalues, oscillation of eigenfunctions, and transformations between such problems are treated. While in the next part of work \cite{Bin1} main objective is to solve inverse problem, so that recover $q$ and other coefficients of problem from given data, which are eigenvalues and norming constants. The norming constants concern "norming"  of  eigenfunctions, which involves a Hilbert space structure that authors set up in work. This allows to view the problem as a standard eigenvalue problem for a selfadjoint operator with compact resolvent and completeness and expansion results are corollaries. For analogous settings see works \cite{Bar,Gul}. We refer here also to \cite{Gok}, where no linearization technique was used, but when $q$is continous (in \cite{Bin,Bin1} $q$ is from $L_1$) more sharpened asymptotic formulas for eigenvalues and eigenfucntions are obtained, and uniform convergence of the continuous functions in series of eigenfunctions is discussed.

But studies of boundary values problems with $\lambda $ included linearly in equation and in boundary conditions or just in boundary conditions for partial diffetrential equations were treated long before than  for scalar Sturm-Liouville problem in \cite{Ful,Wal}. Namely, in 1959 Odhnoff  \cite{Odn} considered spectral problem for Laplace operator in Rimanian manifold when $\lambda $ appears linearly in boundary condition. Later, Ercolano and Schechter \cite{Erc,Erc1} considered $2r$-th order elliptic equation when $\lambda $ included linearly in equation and in boundary conditions. They consider the problem
\begin{equation} \label{eq-1_4}
Ay=\lambda y,{\kern 1pt} {\kern 1pt} {\kern 1pt} {\kern 1pt} on{\kern 1pt} {\kern 1pt} {\kern 1pt} {\kern 1pt} \Omega
\end{equation}

\begin{equation} \label{eq-1_5}
B_{j} y=\lambda C_{j} y,{\kern 1pt} {\kern 1pt} {\kern 1pt} {\kern 1pt} on{\kern 1pt} {\kern 1pt} {\kern 1pt} {\kern 1pt} \Gamma \subset \partial \Omega
\end{equation}
 $A$ is $2r$-th order elliptic operator, $B_j,$ and $C_l$ are linear partial differential operators and $orderB_{j} \ge orderC_{j} $. Moreover, they considered two other, so called domain eigenvalue problem when $\lambda $ appears only in equation,  and boundary eigenvalue problem where $\lambda $  appears only in boundary condition.   They define symmetric operator in $L_{2} (\Omega ){\mathop{\oplus }\nolimits_{}} L_{2} (\Gamma ){\mathop{\oplus }\nolimits_{}} \ldots {\mathop{\oplus }\nolimits_{}} L_{2} (\Gamma )$ and searched its selfadjoint extension and Green function of resolvent. Here we refer also to works by Kozhevnikov, Yakubov \cite{Koz}, Kotko and Krein \cite{Kot}, Kozhevnikov \cite{Koz1}, where considered non-selfadjoint case and asymptotics of eigenvalue distribution and completeness and minimality of system of eigen and associated vectors were studied. In work \cite{Bin21}  Binding, Najman, Langer, Hrinyv consider the case of a formally self-adjoint ellitic operator $A$ of the second order with general than it is in \cite{Odn} and \cite{Erc,Erc1} $\lambda $-linear boundary conditions:
\begin{equation} \label{eq-1_6}
B_{0} u=\lambda B_{1} u{\kern 1pt} {\kern 1pt} {\kern 1pt} {\kern 1pt} on{\kern 1pt} {\kern 1pt} {\kern 1pt} {\kern 1pt} \partial \Omega
\end{equation}
 $B_{0} $ and $B_{1} $ being boundary operators
\[(B_{j} u)(x):=\alpha _{j} (x)\frac{\partial u(x)}{\partial \nu } +\beta _{j} (x)(\gamma u)(x),j=0,1\]
where $\nu (x)$ is the unit normal vector to $\partial \Omega $ pointing out of $\Omega $ at the point $x$ and $\gamma $ is trace operator. When determinant formed by coefficients of boundary conditions
\[\left|\begin{array}{cc} {\alpha _{0} } & {\beta _{0} } \\ {\alpha _{1} } & {\beta _{1} } \\ {} & {} \end{array}\right|=\delta >0\]
 or without loss of generality $\delta =1$ (one arrive at that by dividing both sides of (1.6) by $\sqrt{\delta } $) they introduce Hilbert space
\begin{equation} \label{eq-1_7}
H=L_{2} (\Omega ){\mathop{\oplus }\limits_{}} L_{2} (\partial \Omega )
\end{equation}
 and show symmetricity of operator $A_{0} $ associated with considered problem and defined there. Note that because of $\lambda $ in boundary condition operator corresponding to problem can't be selfadjoint or even symmetric in $L_{2} (\Omega )$. But when writing Green identity for operator defined in \eqref{eq-1_7} the addition of scalar product in second component of direct sum   provides vanishing of scalar products over boundaries giving rise to the relation
\[(A_{0} U,V)-(U,A_{0} V)=0\]
implying symmetricity. If $\delta <0$ then $A_{0} $ is selfadjoint in Krein space defined over $H$ by indefinite metrics
\[[\cdot ,\cdot ]=(\cdot ,\cdot )_{L_{2} (\Omega )} -(\cdot ,\cdot )_{L_{2} (\partial \Omega )} \]
In work \cite{Bin21} spectrum of $A_{0} $ is investigated.

Here we want also a refer to work by Russakovskii \cite{Rus} where boundary conditions involve polynomials on $\lambda $ and by linearizing the problem and introducing  indefinite metrics it was shown that corresponding operator is selfadjoint with respect to that metrics, and operator is selfadjoint in usual sense if and only if ratio of polynomials in $\lambda $ in boundary condtions $R(\lambda )$ maps upper complex plane into itself, which agrees also with results by Walter \cite{Wal} where
\begin{equation} \label{eq-1_8}
R(\lambda )=\frac{a\lambda +b}{c\lambda +d}
\end{equation}
 Posing right definitness condition $\delta =\left|\begin{array}{cc} {a} & {b} \\ {c} & {d} \\ {} & {} \end{array}\right|>0$ in \cite{Asl3} the problem was reduced to spectral problem for selfadjoint operator. But it is well known that $\delta >0$ is sufficient conditon for Mobius transformation $\lambda \to \frac{a\lambda +b}{c\lambda +d} $ to be Herglotz-Nevanlinna type function. Recall that problems with more general Herglotz-Nevanlinna functions of  $\lambda $in boundary condition were treated as indicated above in \cite{Bar,Bin,Bin1,Gul}.

 Note that  there is relation between boundary value problems with $\lambda $ dependent boundary conditions and general theory of selfadjoint extensions with exit from space of minimal and closed operators. It is known that symmetric closed operators always have selfadjoint extensions, (even in case of unequal deficiency numbers) with exit from space. Such problems were subject of studies by Neumark \cite{Neu,Neu1}, Krein \cite{Kre}, in series of works by Shtraus \cite{Sht,Sht1,Sht2,Sht3,Sht4} and others. In \cite{Sht,Sht2} the next  manifold is introduced
\[\tilde{\Omega }_{\lambda } =\left\{\tilde{f}\in D(\tilde{A})/\tilde{A}\tilde{f}-\lambda \tilde{f}\in H\right\}\]
where $\tilde{A}$ is a selfadjoint extension in $\tilde{H}$ of symmetric and closed operator $A$ acting in $H\subset \tilde{H}$. If $P$ projects $\tilde{H}$ to $H$, then in $\Omega _{\lambda } =P\tilde{\Omega }_{\lambda } $ the next operator family depending on $\lambda $
\[B_{\lambda } =\tilde{A}\Omega _{\lambda } \]
is  compression of selfadjoint extension with exit from space to  H and $B_{\lambda } $ is also called quaziselfadjoint extension of $A$.

\noindent Realization of $B_{\lambda } $ to boundary value problems for differential equations is problems with  $\lambda $ dependent boundary conditions. Since $\tilde{A}$ is selfadjoint then for complex $\lambda $ operator $\tilde{A}-\lambda I$($I$-identity operator) has bounded inverse defined everywhere in $\tilde{H}$.Let $R_{\lambda } =P(\tilde{A}-\lambda \tilde{I})^{-1} $. Resolvent of $B_{\lambda } $: $(B_{\lambda } -\lambda I)^{-1} $  called generalized resolvent of $A$   and its selfadjoint extension are related by
\[(B_{\lambda } -\lambda I)^{-1} =R_{\lambda } .\]

 For minimal symmetric differential operator with deficiency indices $(1,1)$ in \cite{Sht4} it was shown that there is one-to-one correspondence between all generalized resolvents of that operatorand Herglotz-Nevanlinna functions. In \cite{Gul1} given description of all selfadjoint extensions with exit from space of minimal closed symmetric operator with deficiency indices $(1,1)$. We refer here also to \cite{Ali,Asl,Asl1,Asl2,AslBay,AslBay1,AslBay2,Asl3,Ryb} where boundary value problems for operator differential equations are considered.

\noindent

\section{Minimal and maximal operators. Selfadjoint extensions.}

\noindent

\noindent Recall that we call  $L^{\prime}_0 $ minimal operator. Show now its symmetricity.

Define in the Hilbert space $L_{2} =L_{2} ((0,1),H)\oplus H^{2} $, the inner product of elements $Y=(y(t),\; y_{1} ,\; y_{2} ),\; Z=(z(t),\; z_{1} ,z_{2} ),\; y(t),\; z(t)\in L_{2} ((0,1),H),\; y_{1} ,\; y_{2} ,\; z_{1} ,\; z_{2} \in H$ by
\begin{equation} \label{eq-2_1}
\left(Y,Z\right)_{H_{1} } =\int _{0}^{1} (y(t),\; z(t))dt+(y_{1} ,z_{1} )+(y_{2} ,z_{2} )
\end{equation}

\begin{lemma}
$L'_{0} $is symmetric operator.
\end{lemma}

\bf{Proof.}
Here and in the future we will take $q\left(t\right)\equiv 0$, since $q\left(t\right)$ is bounded its existence is not essential for  studying domains of minimal, maximal operators and selfadjoint extensions  so one can do it without loss of generality.

 Hence
\[\left(L'_{0} Y,Y\right)=\int _{0}^{1}\left(-y''+Ay,y\right) dt-\left(y\left(1\right),y_{1} \right)+\left(y'\left(1\right),y\left(1\right)\right)-\left(y_{1} ,y\left(1\right)\right)\]
After integrating by parts
\[\left(L'_{0} Y,Y\right)=\left. \left(-y',y\right)\right|_{0}^{1} +\int _{0}^{1}\left(y',y'\right) dt-\left(y\left(1\right),y_{1} \right)+\left(y'\left(1\right),y\left(1\right)\right)-\left(y_{1} ,y\left(1\right)\right)=\]
\[=-\int _{0}^{1}\left(y,y''\right) dt+\left(y\left(1\right),y'\left(1\right)\right)-\left(y\left(1\right),y_{1} \right)-\left(y_{1} ,y\left(1\right)\right)\]
But
\[\left(Y,L'_{0} Y\right)=-\int _{0}^{1}\left(y,y''\right) dt+\left(y\left(1\right),y'\left(1\right)\right)-\left(y\left(1\right),y_{1} \right)-\left(y_{1} ,y\left(1\right)\right)\]
which justifies statement.

 Obviously, $L_0$ is also symmetric as closure of  $L^{\prime}_0.$

 $H_{j} $ denotes completions of $D\left(A^{j} \right)$ $\left(j>0\right)$with respect to scalar product
\[\left(f,g\right)_{j} =\left(A^{j} f,A^{j} g\right)\]
It is known \cite{Tre} that $H_{j} $ is Hilbert space and for $j>k$
\[\left\| \cdot \right\| _{j} \ge \left\| \cdot \right\| _{k} \ge \left\| \cdot \right\|      \left(\left\| \cdot \right\| \equiv \left\| \cdot \right\| _{0} \right)\]
and holds embedding $H_{j} \subseteq H_{k} \subseteq H.$

$H_{-j} $  is space which is adjoint to $H$ with respect to scalar product in $H$. Operator $A$is isometric from $H_{1} $ to $H$. Operator $\tilde{A}$ acting from $H$ to $H_{1} $denotes adjoint of $A.$

 Analogues to \cite{Gor75,Gor1} we can state the next lemma  about the domain of definition of maximal operator. Let
\[\tilde{l}\left(y\right)\equiv -y''+\tilde{A}y+q\left(t\right)y.\]

\begin{lemma}
Domain $D\left(L_{0}^{*} \right)$ of maximal operator $L_{0}^{*} $ consists of elements $Z=\left\{z\left(t\right),z_{1} ,z\left(1\right)\right\}$ where
\begin{equation} \label{eq-2_2}
z\left(t\right)=e^{-\sqrt{\tilde{A}} t} f+e^{-\sqrt{\tilde{A}} \left(1-t\right)} g+\frac{1}{2} \int _{0}^{1}e^{-\sqrt{\tilde{A}} \left|t-s\right|}  A^{-\frac{1}{2} } h\left(s\right)ds,
\end{equation}
\[f\in H_{-\frac{1}{4},}\,\, g\in H_{\frac{1}{2}},\,\,h\left(s\right)\in L_{2} \left(H,\left(0,1\right)\right)\]
and
\[L_{0}^{*} Z=\left\{\tilde{l}z,-z\left(1\right),z^{\prime}\left(1\right)-z_{1} \right\}.\]
Note that domain of definition of $L_{0} $ are those elements from $D\left(L_{0}^{*} \right)$ for which $y\left(1\right)=y^{\prime}\left(1\right)=0.$ If
\[y_{0} =y\left(0\right),y^{\prime}_{0} =y^{\prime}\left(0\right)-\tilde{A}^{\frac{1}{2} } y\left(0\right)\]
it is obvious from \eqref{eq-1_1} that
 $\hat{y}=\tilde{A}^{-\frac{1}{4} } y_{0},$ $\hat{y}^{\prime}=\tilde{A}^{\frac{1}{4} } y^{\prime}_{0} $ are from $H$ and they are regularized boundary values of $y\left(t\right)$ and $y^{\prime}\left(t\right)$.
\end{lemma}

\begin{lemma}
Selfadjoint extensions of $L_{0} $ in $L_{2} $denoted by $L_{s} $ are  given by operations $\left\{\hat{l}y,-y\left(1\right),y'\left(1\right)-y_{1} \right\}$ and boundary conditions at zero
\[\cos C\hat{y}^{\prime}-\sin C\hat{y}=0\]
where C is selfadjoint operator in $H.$
\end{lemma}
When $C=\frac{\pi }{2} I$ we get selfadjoint extension $L'$ defined by
\[y\left(0\right)=0\]
\[L'Y=\left\{\hat{l}y,-y\left(1\right),y_{1} -y\left(1\right)\right\}.\]
Denote by $D\left(L\right)$ set of vectors $\left\{y\left(t\right),y_{1} ,y\left(1\right)\right\}$ where $y\left(t\right)$ satisfies the next conditions:

\begin{enumerate}
\item  $y\left(t\right)$ is twice continuously differentiable in norm of $H$on $\left[0,1\right]$ in strong sense.

\item  $y\left(t\right)\in D\left(A\right)$, $\forall t\in \left[0,1\right]$

\item  $ly\in L_{2} \left(H,\left(0,1\right)\right)$
\end{enumerate}

 Then $D\left(L\right)\subset D\left(L'\right)$and for $Y\in D\left(L\right)$
\[LY=\left\{ly,-y\left(1\right),y_{1} -y\left(1\right)\right\}.\]
Thus $L$ is seladjoint operator in $L_{2} $ with domain consisting of vectors $\left\{y\left(t\right),y_{1} ,y_{1} -y\left(1\right)\right\}$ provided $y\left(t\right)$ satisfies a)-c) and moreover
\[y\left(0\right)=0.\]
Similar to \cite{Ali} the next lemma is true.

\begin{lemma}
Spectrum of selfadjoint extension\textbf{ }$L_{s} $is discrete provided that $A^{-1} \in \sigma _{\infty } $if and only if $\cos C$ is compact operator.
\end{lemma}

Since $\cos \frac{\pi }{2} I=0$(zero operator) is compact, so $L$ is discrete operator.

\section{Eigenvalue asymptotics of operator $\boldsymbol{L}$}

 By our assumption $A$ is selfadjoint, positive definite, and $A^{-1} \in \sigma _{\infty } $.

Thus, its spectrum is discrete, denote them in ascending order $\gamma _{1} \le \gamma _{2} \le \gamma _{3} \le ...$.

Let $\gamma _{k} {\rm \sim }a\; k^{\alpha } ,\; a>0,\; \alpha >0$ and $\varphi _{1} ,\varphi _{2} ,...$ denote the eigenvectors of $A$. From spectral expansion theorem
\[A=\sum _{k=1}^{\infty } \gamma _{k} (\cdot ,\varphi _{k} )\varphi _{k} .\]

Let $L_{0} $ operator obtained from $L$ when $q(t)=0$ so that if $QY=\left\{q(t)y,\; 0,0\right\}$ then $L=L_{0} +Q$.

For coefficients $y_{k} (t)=(y(t),\varphi _{k} )$ in $L_{2} (0,1)$ we have the next spectral problem
\begin{equation} \label{eq-3_1}
-y_{k}^{''} (t)+\gamma _{k} y_{k} (t)=\lambda y_{k} (t),
\end{equation}

\begin{equation} \label{eq-3_2}
y_{k} (0)=0,
\end{equation}

\begin{equation} \label{eq-3_3}
y_{k}^{'} (1)=\left(\lambda -\frac{1}{\lambda } \right)y_{k} (1).
\end{equation}
 Since $L$ and $L_{0} $ both are discrete operators we  may denote their eigenvalues in ascending order as  $\mu _{1} \le \mu _{2} \le ...$ and $\lambda _{1} \le \lambda _{2} \le ...$ respectively. Condition \eqref{eq-3_2} gives that solution of this problem from $L_{2} (0,1)$ is
\[y_{k} (t)=c\sin \sqrt{\lambda -\gamma _{k} } t.\]

For that solution to satisfy condition \eqref{eq-3_3} must be held
\begin{equation} \label{eq-3_4}
\sqrt{\lambda -\gamma _{k} } \cos \sqrt{\lambda -\gamma _{k} } =\left(\lambda -\frac{1}{\lambda } \right)\sin \sqrt{\lambda -\gamma _{k} }
\end{equation}
or
\begin{equation} \label{eq-3_5}
ctg\sqrt{\lambda -\gamma _{k} } =\left(\lambda -\frac{1}{\lambda } \right)\frac{1}{\sqrt{\lambda -\gamma _{k} } }
\end{equation}

Setting $\sqrt{\lambda -\gamma _{k} } =z\; $ gives
\begin{equation} \label{eq-3_6}
z\; ctg\; z=\left(z^{2} +\gamma _{k} -\frac{1}{z^{2} +\gamma _{k} } \right)
\end{equation}

Obviously equation \eqref{eq-3_6} can't have complex roots of form $a+ib,\; a,b\ne 0$, because selfadjoint problem \eqref{eq-3_1}-\eqref{eq-3_3} can't have complex roots. Denoting $z=iy,\; y>0$ gives
\begin{equation} \label{eq-3_7}
ythy=\left(\gamma _{k} -y^{2} -\frac{1}{\gamma _{k} -y^{2} } \right)
\end{equation}
 which has roots $y\sim \sqrt{\gamma _{k} } -\frac{1}{2} $. For  real roots $x_{m,k}$ of \eqref{eq-3_6} we have
\begin{equation} \label{eq-3_8}
\pi m<x_{m,k}<\pi (m+1), m=1,2,..
\end{equation}

We come to the next lemma.

\begin{lemma}
Eigenvalues of $L_{0} $ form two different sequences
\[\lambda _{k,m} =\gamma _{k} +x_{m,k}^{2} ,\; \; x_{m,k} \sim \pi m;\, \, \, \, \, \, \lambda _{k} \sim \sqrt{\gamma _{k} } .\]
\end{lemma}

\begin{lemma}
Eigenvalue distribution function $N(\lambda )$ asymptotically behaves as
\[N(\lambda ){\rm \sim }c\lambda ^{\delta } ,\delta =\left\{\begin{array}{c} {\frac{2\alpha }{\alpha +2} ,\; \; \alpha >2} \\ {\frac{\alpha }{2} ,\; \; \; \; \alpha <2} \\ {1,\; \; \; \; \alpha =2} \end{array}\right. \]
\end{lemma}

Proof might be drown in similar to \cite{Ryb} way and we omit it here.

\begin{lemma}
If $\gamma _{k} {\rm \sim }a\; k^{\alpha } ,\; \alpha >0,\; a>0$, then
\[\lambda _{n} {\rm \sim }\mu _{n} {\rm \sim }dn^{\delta } ,\delta =\left\{\begin{array}{c} {\begin{array}{c} {\frac{\alpha +2}{2\alpha } ,\; \; \; \alpha >2} \\ {\frac{2}{\alpha } ,\; \; \; \; \alpha <2} \end{array}} \\ {1,\; \; \; \; \; \alpha =2} \end{array}\right. \]
\end{lemma}

{Proof} might be derived by the way given in \cite{Ryb}.

\section{Regularized trace of operator $\boldsymbol{L}$ }

Let $\left\{\psi _{n} \right\}$ be orthohormal eigenvectors of operator $L_{0} $. From results in \cite{Mak} and \cite{Sad} we have the next theorem.
\begin{theorem}
If $\alpha >2$, then
\begin{equation} \label{eq-4_1}
{\mathop{\lim }\limits_{m\to \infty }} \sum _{n=1}^{n_{m} } \left(\mu _{n} -\lambda _{n} -\left(Q\psi _{n} ,\psi _{n} \right)\right)=0,
\end{equation}
 where $n_{m} $ some subsequence of natural members defined in \cite{Mak} or \cite{Sad}.
\end{theorem}

First find eigenvectors of $L_{0} $. They are solutions of $L_{0} Y=\lambda Y,$ so eigenvectors  of $L_{0} $  are  of form
\begin{equation} \label{eq-4_2}
C_{m,k} \varphi _{k} \left\{y_{k} (t)-\frac{1}{\lambda _{m,k} } y_{k} (1),\; y_{k} (1)\right\},\; \; k=\overline{1,\infty ,}m=\overline{1,\infty }
\end{equation}
where $\lambda _{m,k} $ are eigenvalues of $L_{0} $ an $y_{k} (t)$ are solutions of \eqref{eq-4_1}-\eqref{eq-4_3}.

To normalize vectors \eqref{eq-4_2} first evaluate norm of $\left\{y_{k} (t)\varphi _{k} ,-\frac{1}{\lambda _{m,k} } y_{k} \eqref{eq-1_1}\varphi _{k} ,y_{k} \eqref{eq-1_1}\varphi _{k} \right\}$ in${\ L}_2$. Hence in virtue of $\left(\varphi _{k} ,\varphi _{k} \right)=1$
\[\left\| \cdot \right\| _{L_{2} }^{2} =\int _{0}^{1} {\mathop{\sin }\nolimits^{2}} x_{m,k} t\; dt+\frac{{\mathop{\sin }\nolimits^{2}} x_{m,k} }{\left(x_{m,k}^{2} +\gamma _{k} \right)^{2} } +\]

\[+{\mathop{\sin }\nolimits^{2}} x_{m,k} =\frac{1}{2} -\int _{0}^{1} \frac{\cos 2x_{m,k} t}{2} dt+\]

\[+\frac{{\mathop{\sin }\nolimits^{2}} x_{m,k} }{\left(x_{m,k}^{2} +\gamma _{k} \right)^{2} } +{\mathop{\sin }\nolimits^{2}} x_{m,k} =\]

\[=\frac{1}{2} -\frac{\sin 2x_{m,k} }{4x_{m,k} } +\frac{{\mathop{\sin }\nolimits^{2}} x_{m,k} }{\left(x_{m,k}^{2} +\gamma _{k} \right)^{2} } +{\mathop{\sin }\nolimits^{2}} x_{m,k} =\]

\[=\frac{x_{m,k} \left(x_{m,k}^{2} +\gamma _{k} \right)^{2} -{\mathop{\sin }\nolimits^{2}} x_{m,k} \left(x_{m,k}^{2} +\gamma _{k} \right)^{2} /2}{2x_{m,k} \left(x_{m,k}^{2} +\gamma _{k} \right)^{2} } +\]

\begin{equation} \label{eq-4_3}
+\frac{2x_{m,k} {\mathop{\sin }\nolimits^{2}} x_{m,k} +2x_{m,k} \left(x_{m,k}^{2} +\gamma _{k} \right)^{2} {\mathop{\sin }\nolimits^{2}} x_{m,k} }{2x_{m,k} \left(x_{m,k}^{2} +\gamma _{k} \right)^{2} }
\end{equation}
Taking into consideration replacement $\sqrt{\lambda -{\gamma }_k}=x_{m.k}$ in \eqref{eq-3_4} gives
\begin{equation} \label{eq-4_4}
\sin x_{m,k} \left(x_{m,k}^{2} +\gamma _{k} \right)=x_{m,k} \cos x_{m,k} +\frac{\sin x_{m,k} }{x_{m,k}^{2} +\gamma _{k} }
\end{equation}
and  using  this expression  in the last term from \eqref{eq-4_3}  after simplifications we have
\[\left\| \cdot \right\| ^{2} =\frac{x_{m,k} \left(x_{m,k}^{2} +\gamma _{k} \right)^{2} -\left(x_{m,k}^{2} +\gamma _{k} \right)^{2} \frac{{\mathop{\sin }\nolimits^{2}} x_{m,k} }{2} }{2x_{m,k} \left(x_{m,k}^{2} +\gamma _{k} \right)^{2} } +\]

\begin{equation} \label{eq-4_5}
+\frac{4x_{m,k} {\mathop{\sin }\nolimits^{2}} x_{m,k} +2x_{m,k}^{2} \sin x_{m,k} \cos x_{m,k} \left(x_{m,k}^{2} +\gamma _{k} \right)}{2x_{m,k} \left(x_{m,k}^{2} +\gamma _{k} \right)^{2} } .
\end{equation}

Denote numerator of \eqref{eq-4_5} by $H_{m,k} .$ So normed eigenvectors are of the form
\[\psi _{m,k} =\frac{\sqrt{2x_{m,k} } \left(x_{m,k}^{2} +\gamma _{k} \right)}{\sqrt{H_{m,k} } } \times \]

\begin{equation} \label{eq-4_6}
\times \left\{\sin (x_{m,k} \; t)\varphi _{k} ,-\frac{1}{x_{m,k}^{2} +\gamma _{k} } \sin (x_{m,k} )\varphi _{k} ,\; \sin x_{m,k} \varphi _{k} \right\}
\end{equation}
 Prove the next lemma which will show that value of limit in \eqref{eq-4_1} doesn't depend on choice of $n_{m} .$

\begin{lemma}
If operator-valued function $q(t)$ together with conditions from section 1 satisfies also the next ones:

1) $q(t)$ has the second selfadjoint weak derivatives till order two $\left(q^{(l)} (t)=\left[q^{(l)} (t)\right]^{*} ,\; l=0,1,2\right)\; $ in $H$;
\[2) \sum _{k=1}^{\infty } \left|\left(q^{(l)} (t)\varphi _{k} ,\; \varphi _{k} \right)\right|<const,\; \; \forall t\in 0,1];\]
\[3) \int _{0}^{1} \left(q(t)f,f\right)dt=0,\; \forall f\in H,\]

then
\[\sum _{k=1}^{\infty } \sum _{m=1}^{\infty } \left|\left(Q\psi _{k,m} ,\psi _{k,m} \right)\right|=\]

\[=\sum _{k=1}^{\infty } \sum _{m=1}^{\infty } \left|\int _{0}^{1} \frac{x_{m,k} \left(x_{m,k}^{2} +\gamma _{k} \right)^{2} {\mathop{\sin }\nolimits^{2}} x_{m,k} t}{H_{m,k} } q_{k} (t)dt\right|<\infty \]
where $q_{k} (t)=(q(t)\varphi _{k} ,\; \varphi _{k} ).$
\end{lemma}
\bf{Proof.}
By condition 3) it is suffices to perform  twice integration by parts  the next integral
\[\int _{0}^{1} \cos \left(2x_{m,k} t\right)q_{k} (t)dt=\frac{\sin (2x_{m,k} )}{2x_{m,k} } q_{k} (1)-\]

\[-\int _{0}^{1} \frac{\cos (2x_{m,k} t)}{(2x_{m,k} )^{2} } q_{k}^{''} (t)dt+\frac{\cos 2x_{m,k} }{(2x_{m,k} )^{2} } q_{k}^{'} (1)-\frac{q_{k}^{'} (0)}{(2x_{m,k} )^{2} } \]
Here we also use condition 1). The last equality together with conditions 2), 3) and asymptotics $x_{m,k} $ proves validity of statement.

Now we may state that
\begin{equation} \label{eq-4_7}
\sum _{n=1}^{n_{m} } \left(\mu _{n} -\lambda _{n} \right)
\end{equation}
converges as $m\to \infty ,$ its value doesn't depend on choice of $\left\{n_{m} \right\}$ and
\[{\mathop{\lim }\limits_{m\to \infty }} \sum _{n=1}^{n_{m} } \left(\mu _{n} -\lambda _{n} \right)=\sum _{k=1}^{\infty } \sum _{m=1}^{\infty } \left(Q\psi _{m,k} ,\psi _{m,k} \right)_{L_{2} } =\]

\begin{equation} \label{eq-4_8}
=-\sum _{k=1}^{\infty } \sum _{m=1}^{\infty } \int _{0}^{1} \frac{\left(x_{m,k}^{2} +\gamma _{k} \right)^{2} \cos (2x_{m,k} t)q_{k} (t)}{H_{m,k} } dt
\end{equation}
Call $\lim\limits_{m\to \infty } \sum\limits_{n=1}^{n_{m} }$,$\left(\mu_{n} -\lambda_{n} \right)\; $ the first regularized trace of $L$ and evaluate it.

 For that reason prove the next theorem.
\begin{theorem}
If $\alpha >2$ and hold conditions 1)-3) then
\begin{equation} \label{eq-4_9}
{\mathop{\lim }\limits_{m\to \infty }} \sum _{n=1}^{n_{m} } \left(\mu _{n} -\lambda _{n} \right)=-\frac{tr\; q(0)+tr\; q(1)}{4} .
\end{equation}
\end{theorem}
\bf{Proof.}
For evaluating the value of  inner series in \eqref{eq-4_8} we have to investigate the asymptotic behavior of sum
\[M_{N} (t)=\sum _{m=1}^{N} \frac{x_{m,k} \left(x_{m,k}^{2} +\gamma _{k} \right)^{2} \cos (2x_{m,k} t)q_{k} (t)}{M_{m,k} } \]
 as $N\to \infty $. For that purpose we will find poles are $g(x)$ whose poles are $x_{m,k} \; $and residues at that poles yield the terms of sum \eqref{eq-4_7}. Then integrating that function along extending contours containing $x_{m,k} $ and using asymptotics of $g(z)$ we will get ${\mathop{\lim }\limits_{N\to \infty }} M_{N} (t)$.

Consider the function
\[g(z)=\frac{z\cos (2\; zt)}{\left(z\; ctg-z^{2} -\gamma _{k} +\frac{1}{z^{2} +\gamma _{k} } \right){\mathop{\sin }\nolimits^{2}} z} \]
 for each fixed $t$. As it is seen from \eqref{eq-4_5} its simple poles are $x_{m,k} $ and $\pi m$.

Evaluate its residues at that poles.

First evaluate $f^{'} \left(x_{m,k} \right)\; $where $f(z)=\left(z\; ctg\; z-z^{2} -\gamma _{k} +\frac{1}{z^{2} +\gamma _{k} } \right)$
\[f^{'} (x_{m,k} )=ctg\; x_{m,k} -\frac{x_{m,k} }{{\mathop{\sin }\nolimits^{2}} x_{m,k} } -2x_{m,k} -\]

\[-\frac{2x_{m,k} }{\left(x_{m,k}^{2} +\gamma _{k} \right)^{2} } =\frac{\sin x_{m,k} \cos x_{m,k} \left(x_{m,k}^{2} +\gamma _{k} \right)^{2} -x_{m,k} \left(x_{m,k}^{2} +\gamma _{k} \right)^{2} }{{\mathop{\sin }\nolimits^{2}} (x_{m,k} )\left(x_{m,k}^{2} +\gamma _{k} \right)^{2} } +\]

\[+\frac{-2x_{m,k} \left(x_{m,k}^{2} +\gamma _{k} \right)^{2} {\mathop{\sin }\nolimits^{2}} x_{m,k} -2x_{m,k} {\mathop{\sin }\nolimits^{2}} x_{m,k} }{{\mathop{\sin }\nolimits^{2}} (x_{m,k} )\left(x_{m,k}^{2} +\gamma _{k} \right)^{2} } =\]

\[=\frac{\frac{\sin 2x_{m,k} }{2} \left(x_{m,k}^{2} +\gamma _{k} \right)^{2} -x_{m,k} \left(x_{m,k}^{2} +\gamma _{k} \right)^{2} }{{\mathop{\sin }\nolimits^{2}} (x_{m,k} )\left(x_{m,k}^{2} +\gamma _{k} \right)^{2} } +\]

\[+\frac{-2x_{m,k} \sin x_{m,k} \left(x_{m,k}^{2} +\gamma _{k} \right)\left(x_{m,k} \cos x_{m,k} +\frac{\sin x_{m,k} }{x_{m,k}^{2} +\gamma _{k} } \right)-2x_{m,k} {\mathop{\sin }\nolimits^{2}} x_{m,k} }{{\mathop{\sin }\nolimits^{2}} (x_{m,k} )\left(x_{m,k}^{2} +\gamma _{k} \right)^{2} } =\]
\[=\frac{\frac{\sin 2x_{m,k} }{2} \left(x_{m,k}^{2} +\gamma _{k} \right)^{2} -x_{m,k} \left(x_{m,k}^{2} +\gamma _{k} \right)^{2} -2x_{m,k} \left(x_{m,k}^{2} +\gamma _{k} \right)^{2} \sin x_{m,k} \cos x_{m,k} -4x_{m,k} \sin ^{2} x_{m,k} }{\sin ^{2} x_{m,k} \left(x_{m,k}^{2} +\gamma _{k} \right)^{2} } \]
Here we have used  \eqref{eq-4_4}. Thus  we get

\[{resg(z)}_{z=x_{m.k}=}\frac{x_{m,k} \cos (2x_{m,k} t)}{f^{'} (x_{m,k} ){\mathop{\sin }\nolimits^{2}} x_{m,k} } \]

\noindent which are the terms of sum \eqref{eq-4_7}. For points $\pi m$ we have
\[{\mathop{res}\limits_{z=\pi m}} g(z)=\frac{\pi _{m} \cos (2\pi mt)}{\cos \pi m\left(\pi m\cos \pi m-(\pi m\right)^{2} \sin \pi m-\gamma _{k} \sin \pi m-\frac{\sin \pi m}{\left(\pi m\right)^{2} +\gamma _{k} } )} \]

\[=\cos 2\pi mt.\]

Denote by $S_{N} (t)$ the next sum
\[S_{N} (t)=\sum _{m=1}^{N} \cos 2\pi mt\]

We will apply about residues. Take as contour of integration rectangle with vertices at points $\pm iB,\; A_{N} \pm iB$, further $B$ will tend to infinity and $A_{N} =\pi n+\frac{\pi }{2} $. For such values of $A_{N} $ from asymptotics of $x_{m,k} $ for big $m$ values
\[x_{N,k} <A_{N} <x_{N+1,k} \]
 and
\[\pi N<A_{N} <\pi (N+1).\]

 Let contour by pass origin on the left along small semicirce with radius $r$. Integrating $g(z)$ along that contour, by theorem about residues, we will have
\[M_{N} (t)+S_{N} (t)=\frac{1}{2\pi i} \left[\int _{A_{N} -iB}^{A_{N} +iB} g(z)dz+{\mathop{\lim }\limits_{r\to 0}} \int _{{\mathop{\left|z\right|=r}\limits_{-\frac{\pi }{2} <\varphi <\frac{\pi }{2} }} } g(z)dz+\int _{iB}^{A_{N} +iB} g(z)dz+\int _{iB}^{A_{N} -iB} g(z)dz\right].\]

If $z={\mathop{u+tv}\limits_{u\ge 0}} $, then order of $g(z)$ is $O$ $\left(\frac{1}{e^{2\left|v\right|(1-t)} } \right),$ thus integrals along upper and lower sides of integration tends zero as $B\to \infty $. It might be easily  shown that integral along small semicircle vanishes when $r\to 0$, because $g(z){\rm \sim }re^{i\varphi } t^{2} \to 0\; \; \left(t\in \left[0,1\right]\right)$

For $\left|z\right|\to \infty $ from asymptotics of $g(z)$ we can show that the first integral vanishes as $B\to \infty $. Hence,
\begin{equation} \label{eq-4_10}
{\mathop{\lim }\limits_{N\to \infty }} \int _{0}^{1} M_{N} (t)q_{k} (t)dt=-{\mathop{\lim }\limits_{N\to \infty }} \int _{0}^{1} S_{N} (t)q_{k} (t)dt.
\end{equation}

But
\begin{equation} \label{eq-4_11}
{\mathop{\lim }\limits_{N\to \infty }} \int _{0}^{1} S_{N} (t)q_{k} (t)dt=-\frac{q_{k} (1)+q_{k} (0)}{4} .
\end{equation}

Now from \eqref{eq-4_8} and \eqref{eq-4_10}, \eqref{eq-4_11} summing for $k$ from $1$ to $\infty $ we get validity of theorem.

\end{document}